\documentclass{amsart}
\usepackage{amssymb}
\usepackage{amsmath}
\usepackage{amsfonts}

\newtheorem{theorem}{Theorem}
\theoremstyle{plain}

\numberwithin{equation}{section}
\input{tcilatex}

\begin{document}
\title[Intrinsic characterization of locally C*- and locally JB-Algebras]{On
Intrinsic Characterization of Real Locally C*- and Locally JB-Algebras}
\author{Alexander A. Katz}
\address{Dr. Alexander A. Katz, Department of Mathematics and Computer
Science, St. John's College of Liberal Arts and Sciences, St. John's
University, 300 Howard Avenue, DaSilva Academic Center 314, Staten Island,
NY 10301, USA}
\email{katza@stjohns.edu}
\author{Oleg Friedman}
\address{Oleg Friedman, Department of Mathematical Sciences, University of
South Africa, P.O. Box 392, Pretoria 0003, South Africa}
\email{friedman001@yahoo.com}
\curraddr{Oleg Friedman, Department of Mathematics and Computer Science, St.
John's College of Liberal Arts and Sciences, St. John's University, 8000
Utopia Parkway, St. John's Hall 334, Queens, NY 11439, USA}
\email{fridmano@stjohns.edu}
\thanks{The second author is thankful to the first author and Dr. Louis E.
Labuschagne (University of South Africa, Pretoria, South Africa) for
constant help, support and discussions.}
\date{October 13-14, 2007}
\dedicatory{Dedicated to the memory of Professor George Bachman, Polytechnic
University, Brooklyn, NY, USA.}

\begin{abstract}
In the present paper we obtain an intrinsic characterization of real locally
C*-algebras (projective limits of projective families of real C*-algebras)
among complete real lmc *-algebras, and of locally JB-algebras (projective
limits of projective families of JB-algebras) among complete fine Jordan
locally multiplicatively-convex topological algebras.
\end{abstract}

\maketitle

\section{Introduction}

Banach associative regular *-algebras over $%
%TCIMACRO{\U{2102} }%
%BeginExpansion
\mathbb{C}
%EndExpansion
$, so called \textit{C*-algebras}, were first introduces in 1940's by
Gelfand and Naimark in the paper \cite{GelfandNaimark43}. Since then these
algebras were studied extensively by various authors, and now, the theory of
C*-algebras is a huge part of Functional Analysis which found applications
in almost all branches of Modern Mathematics and Theoretical Physics. For
the basics of the theory of C*-algebras, see for example Pedersen's
monograph \cite{Pedersen79}.

The real analogues of complex C*-algebras, so called \textit{real C*-algebras%
}, which are real Banach *-algebras with regular norms such that their
complexifications are complex\textit{\ C*-algebras, }were studied in
parallel by many authors. For the current state of the basic theory of real
C*-algebras, see Li's monograph \cite{Li03}.

The real Jordan analogues of complex C*-algebras, so called \textit{%
JB-algebras}, were first defined by Alfsen, Schultz and St\o rmer in \cite%
{AlfsenSchultzStoermer78} as the real Banach--Jordan algebras satisfying for
all pairs of elements $x$ and $y$ the inequality of fineness%
\begin{equation*}
\left\Vert x^{2}+y^{2}\right\Vert \geq \left\Vert x\right\Vert ^{2},
\end{equation*}%
and regularity identity 
\begin{equation*}
\left\Vert x^{2}\right\Vert =\left\Vert x\right\Vert ^{2}.
\end{equation*}%
The basic theory of JB-algebras is fully treated in monograph of
Hanche-Olsen and St\o rmer \cite{Hanche-OlsenStoermer84}. If $A$ is a
C*--algebra, or a real C*-algebra, then the self-adjoint part $A_{sa}$ of $A$
is a JB-algebra under the Jordan product 
\begin{equation*}
x\circ y=\frac{(xy+yx)}{2}.
\end{equation*}%
Closed subalgebras of $A_{sa}$, for some C*-algebra or real C*-algebra $A$,
become relevant examples of JB-algebras, and are called \textit{JC-algebras}.

Complete locally multiplicatively-convex algebras or equivalently, due to
Arens-Michael Theorem, projective limits of projective families of Banach
algebras, were first studied by Arens in \cite{Arens52} and Michael in \cite%
{Michael52}. They were since studied by many authors under different names.
In particular, projective limits of projective families of C*-algebras were
studied by Inoue in \cite{Inoue71}, Apostol in \cite{Apostol71}, Schm\"{u}%
dgen in\textbf{\ }\cite{Schmuedgen75}, Phillips in \cite{Phillips88}, Bhatt
and Karia in \cite{BhattKaria93}, etc. We will follow Inoue \cite{Inoue71}
in the usage of the name \textit{locally C*-algebras} for these topological
algebras. The current state of the basic theory of locally C*-algebras is
treated in the monograph of Fragoulopoulou\textbf{\ \cite{Fragoulopoulou05}.}

Topological algebras which are projective limits of projective families of
real C*-algebras under the name of \textit{real locally C*-algebras}, and
projective limits of projective families of JB-algebras under the name of 
\textit{locally JB-algebras} were first introduced by Katz and Friedman in 
\cite{KatzFriedman06}.

Bhatt and Karia in \cite{BhattKaria93} studied the structure of locally
C*-algebras. They obtain the following characterization of locally
C*-algebras among complete lmc *-algebras (topological algebras which are
projective limits of projective families of complex Banach *-algebras):

\begin{theorem}[Bhatt and Karia \protect\cite{BhattKaria93}]
Let $\mathfrak{A}$ be a complex complete lmc *-algebra. Then $\mathfrak{A}$
is a locally C*-algebra iff $\mathfrak{A}$ contains a *-subalgebra $%
\mathfrak{B}$ such that:

1. $\mathfrak{B}$ is a C*-algebra with some norm $\left\Vert .\right\Vert _{%
\mathfrak{B}};$

\ \ \ \ \ \ \ \ \ \ \ \ \ \ \ \ \ \ \ \ \ \ \ \ \ \ \ \ \ \ \ \ \ \ \ \ \ \
\ \ \ \ and

2. the inclusion 
\begin{equation*}
(\mathfrak{B},\left\Vert .\right\Vert _{\mathfrak{B}})\rightarrow \mathfrak{A%
},
\end{equation*}%
is a continuous embedding with dense range.

Further, if the unit ball of $\mathfrak{B},$ 
\begin{equation*}
\mathfrak{B}_{\mathbf{1}}=\{x\in \mathfrak{B}:\left\Vert x\right\Vert _{%
\mathfrak{B}}\leq 1\},
\end{equation*}%
is closed in $\mathfrak{A}$ in projective topology of $\mathfrak{A}$, then 
\begin{equation*}
\mathfrak{B}=\mathfrak{A}_{b},
\end{equation*}%
where by $\mathfrak{A}_{b}$ we mean the bounded part of $\mathfrak{A}$ (see
below Section 2 for precise definitions). \ \ $\square $
\end{theorem}

\begin{proof}
The \textquotedblleft only if\textquotedblright\ part is due to Apostol (see 
\cite{Apostol71}), and, by different methods, to Schm\"{u}dgen (see \cite%
{Schmuedgen75}).and Phillips (see \cite{Phillips88}). The \textquotedblleft
if\textquotedblright\ part is based on numerical range theory in lmc
*-algebras, developed by Giles and Koehler in \cite{GilesKoehler73}.
\end{proof}

The present paper is devoted to the presentation of analogues of Theorem 1
for real locally C*-algebras and locally JB-algebras. In particular, we give
an intrinsic characterization of a real locally C*-algebra as if and only if
it is a complete real lmc *-algebra (projective limit of a projective family
of real Banach *-algebras) with a continuously embedded dense *-subalgebra
which is a real C*-algebra under some norm, as well as an intrinsic
characterization of a locally JB-algebra as if and only if it is a complete
fine locally Banach-Jordan algebra (projective limit of a projective family
of fine Banach-Jordan algebras) with a continuously embedded dense Jordan
subalgebra which is a JB-algebra under some norm.

\section{Preliminaries}

Let us briefly recall some of the basic material from the aforementioned
sources one needs to comprehend what follows.

A Hausdorff topological vector space over the field of $%
%TCIMACRO{\U{211d} }%
%BeginExpansion
\mathbb{R}
%EndExpansion
$ or $%
%TCIMACRO{\U{2102} }%
%BeginExpansion
\mathbb{C}
%EndExpansion
$, in which any neighborhood of the zero element contains a convex
neighborhood of the zero element; in other words, a topological vector space
is a \textit{locally convex space} if and only if the topology of is a
Hausdorff locally convex topology.

A number of general properties of locally convex spaces follows immediately
from the corresponding properties of locally convex topologies; in
particular, subspaces and Hausdorff quotient spaces of a locally convex
space, and also products of families of locally convex spaces, are
themselves locally convex spaces. Let $\Lambda $ be an upward directed set
of indices and a family 
\begin{equation*}
\{E_{\alpha },\alpha \in \Lambda \},
\end{equation*}%
of locally convex spaces (over the same field) with topologies 
\begin{equation*}
\{\tau _{\alpha },\alpha \in \Lambda \}.
\end{equation*}%
Suppose that for any pair $(\alpha ,\beta )$, 
\begin{equation*}
\alpha \leq \beta ,
\end{equation*}%
$\alpha ,\beta \in \Lambda $, there is defined a continuous linear mapping 
\begin{equation*}
g_{\alpha }^{\beta }:E_{\beta }\rightarrow E_{\alpha }.
\end{equation*}

A family 
\begin{equation*}
\{E_{\alpha },\alpha \in \Lambda \}
\end{equation*}%
is called \textit{projective}, if for each triplet $(\alpha ,\beta ,\gamma
), $ 
\begin{equation*}
\alpha \leq \beta \leq \gamma ,
\end{equation*}%
$\alpha ,\beta ,\gamma \in \Lambda ,$ 
\begin{equation*}
g_{\alpha }^{\gamma }=g_{\beta }^{\gamma }\circ g_{\alpha }^{\beta },
\end{equation*}%
and for each $\alpha \in \Lambda ,$ 
\begin{equation*}
g_{\alpha }^{\alpha }=Id.
\end{equation*}

Let $E$\ be the subspace of the product 
\begin{equation*}
\dprod\limits_{\alpha \in \Lambda }E_{\alpha },
\end{equation*}%
whose elements 
\begin{equation*}
x=(x_{\alpha }),
\end{equation*}%
satisfy the relations 
\begin{equation*}
x_{\alpha }=g_{\alpha }^{\beta }(x_{\beta }),
\end{equation*}%
for all $\alpha \leq \beta $. The space $E$\ is called the \textit{%
projective limit} of the projective family $E_{\alpha },$ $\alpha \in
\Lambda ,$\ with respect to the family $(g_{\alpha }^{\beta }),$ $\alpha
,\beta \in \Lambda $\ and is denoted by 
\begin{equation*}
\lim g_{\alpha }^{\beta }E_{\beta },
\end{equation*}%
or 
\begin{equation*}
\underset{\longleftarrow }{\lim }E_{\alpha }.
\end{equation*}%
The topology of $E$\ is the \textit{projective topology} with respect to the
family 
\begin{equation*}
(E_{\alpha },\tau _{\alpha },\pi _{\alpha }),
\end{equation*}%
$\alpha \in \Lambda $, where $\pi _{\alpha },$ $\alpha \in \Lambda ,$\ is
the restriction to the subspace $E$\ of the projection 
\begin{equation*}
\widehat{\pi }_{\alpha }:\dprod\limits_{\beta \in \Lambda }E_{\beta
}\rightarrow E_{\alpha },
\end{equation*}%
and 
\begin{equation*}
\pi _{\beta }=g_{\alpha }^{\beta }\circ \pi _{\alpha },
\end{equation*}%
$\forall \alpha ,\beta \in \Lambda ,$

When you take instead of $E_{\alpha },$ $\alpha \in \Lambda ,$ a projective
family of algebras, *-algebras, Jordan algebras, etc., you naturally get a
correspondent algebra, *-algebra or Jordan algebra structure in the
projective limit algebra 
\begin{equation*}
E=\underset{\longleftarrow }{\lim }E_{\alpha }.
\end{equation*}

Let $E$ be a vector space. A real function $p:E\rightarrow 
%TCIMACRO{\U{211d} }%
%BeginExpansion
\mathbb{R}
%EndExpansion
$ on $E$ is called a \textit{seminorm}, if:

\textit{1).} $p(x)\geq 0,$ $\forall x\in E;$

\textit{2).} $p(\lambda x)=\left\vert \lambda \right\vert p(x),$ $\forall
\lambda \in 
%TCIMACRO{\U{211d} }%
%BeginExpansion
\mathbb{R}
%EndExpansion
$ or $%
%TCIMACRO{\U{2102} }%
%BeginExpansion
\mathbb{C}
%EndExpansion
$, and $x\in E;$

\textit{3).} $p(x+y)\leq p(x)+p(y),$ $\forall x,y\in E.$

One can see that 
\begin{equation*}
p(\mathbf{0})=0.
\end{equation*}%
If 
\begin{equation*}
p(x)=0,
\end{equation*}%
implies 
\begin{equation*}
x=\mathbf{0},
\end{equation*}%
seminorm is called a \textit{norm and is usually denoted by }$\left\Vert
.\right\Vert $. If a space with a norm is complete, it is called a \textit{%
Banach space}.

Let $(E,p)$ be a seminormed space, and 
\begin{equation*}
N_{p}=\ker (p)=p^{-1}\{0\}.
\end{equation*}%
The quotient space $E/N_{p}$ is a linear space and the function 
\begin{equation*}
\left\Vert .\right\Vert _{p}:E/N_{p}\rightarrow 
%TCIMACRO{\U{211d} }%
%BeginExpansion
\mathbb{R}
%EndExpansion
_{+}:
\end{equation*}%
\begin{equation*}
x_{p}=x+N_{p}\rightarrow \left\Vert x_{p}\right\Vert _{p}=p(x),
\end{equation*}%
is a well defined norm on $E/N_{p}$ induced by the seminorm $p.$ The
corresponding quotient normed space will be denoted by $E/N_{p},$ and the
Banach space completion of $E/N_{p}$ by $E_{p}.$ One can easily see that $%
E_{p}$ is the Hausdorff completion of the seminormed space $(E,p).$

The algebras considered below will be without the loss of generality unital.
If the algebra does not have an identity, it can be adjoint by a usual
unitialization procedure.

A Jordan algebra is an algebra $E$ in which the identities 
\begin{equation*}
x\circ y=y\circ x,
\end{equation*}%
\begin{equation*}
x^{2}\circ (y\circ x)=(x^{2}\circ y)\circ x,
\end{equation*}%
hold.

If $E$ is an algebra, the seminorm $p$ on $E$ compatible with the
multiplication of $E$, in the sense that 
\begin{equation*}
p(xy)\leq p(x)p(y),
\end{equation*}%
$\forall x,y\in E,$ is called \textit{submultiplicative }or \textit{%
m-seminorm}.

For submultiplicative seminorm on a Jordan algebra $E$, the following
inequality holds: 
\begin{equation*}
p(x\circ y)\leq p(x)p(y),
\end{equation*}%
$\forall x,y\in E.$ A seminorm on a Jordan algebra $E$\ is called \textit{%
fine}, if the following inequality holds: 
\begin{equation*}
p(x^{2}+y^{2})\geq p(x^{2}),
\end{equation*}%
$\forall x,y\in E.$

A \textit{Banach-Jordan algebra} is Jordan algebra which is as well a Banach
algebra.

Let $E$ be an algebra. A subset $U$ of $E$ is called \textit{multiplicative}
or \textit{idempotent}, if 
\begin{equation*}
UU\subseteq U,
\end{equation*}%
in the sense that $\forall x,y\in U,$ the product 
\begin{equation*}
xy\in U.
\end{equation*}

If $p$ is an m-seminorm on $E$ the unit semiball $U_{p}(1)$ corresponding to 
$p$, that is 
\begin{equation*}
U_{p}(1)=\{x\in E:p(x)\leq 1\},
\end{equation*}%
and one can see that this set is multiplicative. Moreover, $U_{p}(1)$ is an
absolutely-convex (balanced and convex).absorbing subset of $E.$ It is known
that given an absorbing absolutely-convex subset 
\begin{equation*}
U\subset E,
\end{equation*}%
the function%
\begin{equation*}
p_{U}:E\rightarrow 
%TCIMACRO{\U{211d} }%
%BeginExpansion
\mathbb{R}
%EndExpansion
_{+}:
\end{equation*}%
\begin{equation*}
x\rightarrow p_{U}(x)=\inf \{\lambda >0:x\in \lambda U\},
\end{equation*}%
called \textit{gauge} or \textit{Minkowski functional} of $U,$ is a
seminorm. One can see that a real-valued function $p$ on the algebra $E$ is
an m-seminorm iff 
\begin{equation*}
p=p_{U},
\end{equation*}%
for some absorbing, absolutely-convex and multiplicative subset 
\begin{equation*}
U\subset E.
\end{equation*}%
In fact, one can take 
\begin{equation*}
U=U_{p}(1).
\end{equation*}

By \textit{topological algebra} we mean a topological vector space which is
also an algebra, such that the ring multiplication is separately continuous.
A topological algebra $E$ is often denoted by $(E,\tau ),$ where $\tau $ is
the topology of the underlying topological vector space of $E.$ The topology 
$\tau $ is determined by a \textit{fundamental }$0$\textit{-neighborhood
system}, say $\mathcal{B}$, consisting of absorbing, balanced sets with the
property 
\begin{equation*}
\forall V\in \mathcal{B}\exists U\in \mathcal{B},
\end{equation*}%
satisfying the condition $U+U\subseteq V.$ Since translations by $y$ in $%
(E,\tau )$, i.e. the maps 
\begin{equation*}
x\rightarrow x+y:
\end{equation*}%
\begin{equation*}
(E,\tau )\rightarrow (E,\tau ),
\end{equation*}%
$y\in E,$ are homomorphisms, an $x$-neighborhood in $(E,\tau )$ is of the
form 
\begin{equation*}
x+V,
\end{equation*}%
with $V\in \mathcal{B}$. A closed, absorbing and absolutely-convex subset of
a topological algebra $(E,\tau )$ is called \textit{barrel}. An \textit{%
m-barrel} is a multiplicative barrel of $(E,\tau ).$

A \textit{locally convex algebra} is a topological algebra in which the
underlying topological vector space is a locally convex space. The topology $%
\tau $\ of a locally convex algebra $(E,\tau )$ is defined by a fundamental $%
0$-neighborhood system consisting of closed absolutely-convex sets.
Equivalently, the same topology $\tau $\ is determined by a family of
nonzero seminorms. Such a family, say 
\begin{equation*}
\Gamma =\{p\},
\end{equation*}%
or, for distinction purposes 
\begin{equation*}
\Gamma _{E}=\{p\},
\end{equation*}%
is always assumed without a loss of generality \textit{saturated}. That is,
for any finite subset 
\begin{equation*}
F\subset \Gamma ,
\end{equation*}%
the seminorm 
\begin{equation*}
p_{F}(x)=\underset{p\in F}{\max }p(x),
\end{equation*}%
$x\in E,$ again belongs to $\Gamma .$ Saying that 
\begin{equation*}
\Gamma =\{p\},
\end{equation*}%
is a \textit{defining family of seminorms} for a locally convex algebra $%
(E,\tau ),$ we mean that $\Gamma $ is a saturated family of seminorms
defining the topology $\tau $ on $E.$ That is 
\begin{equation*}
\tau =\tau _{\Gamma },
\end{equation*}%
with $\tau _{\Gamma }$ completely determined by a fundamental $0$%
-neighborhood system given by the $\varepsilon $-semiballs 
\begin{equation*}
U_{p}(\varepsilon )=\varepsilon U_{p}(\varepsilon )=\{x\in E:p(x)\leq
\varepsilon \},
\end{equation*}%
$\varepsilon >0,$ $p\in \Gamma $. More precisely, for each $0$-neighborhood 
\begin{equation*}
V\subset (E,\tau ),
\end{equation*}%
there is an $\varepsilon $-semiball $U_{p}(\varepsilon )$, $\varepsilon >0,$ 
$p\in \Gamma ,$ such that 
\begin{equation*}
U_{p}(\varepsilon )\subseteq V.
\end{equation*}%
The neighborhoods $U_{p}(\varepsilon )$, $\varepsilon >0,$ $p\in \Gamma ,$
are called \textit{basic }$0$\textit{-neighborhoods}.

A locally C*-algebra (real locally C*-algebra, resp. locally JB-algebra) is
a projective limit of projective family of C*-algebras (real C*-algebras,
resp. JB-algebras). This is equivalent for locally C*- and real locally
C*-algebras to the requirement that the family of defining continuous
seminorms be regular: 
\begin{equation*}
p(x^{\ast }x)=p(x)^{2}.
\end{equation*}
In the case of locally JB-algebras this is equivalent to the requirement
that the family of defining continuous seminorms be fine and regular:

\begin{equation*}
p(x^{2}+y^{2})\geq p(x^{2}),
\end{equation*}
and

\begin{equation*}
p(x^{2})=p(x)^{2},
\end{equation*}%
$\forall p\in \Gamma ,$ $x,y\in E.$

For a locally C*-algebra (real locally C*-algebra, resp. locally JB-algebra) 
$E$, by the bounded part we mean the subalgebra 
\begin{equation*}
E_{b}=\{x\in E:\left\Vert x\right\Vert _{\infty }=\underset{p\in \Gamma (E)}{%
\sup }p(x)<\infty \}.
\end{equation*}

\section{Intrinsic characterization of real locally C*-algebras}

In the current section we present a real analogue of Theorem 1.

\begin{theorem}
Let $A$ be a complex complete real lmc *-algebra. Then $A$ is a real locally
C*-algebra iff $A$ contains a *-subalgebra $B$ such that:

1. $B$ is a real C*-algebra with some norm $\left\Vert .\right\Vert _{B};$

\ \ \ \ \ \ \ \ \ \ \ \ \ \ \ \ \ \ \ \ \ \ \ \ \ \ \ \ \ \ \ \ \ \ \ \ \ \
\ \ \ \ and

2. the inclusion 
\begin{equation*}
(B,\left\Vert .\right\Vert _{B})\rightarrow A,
\end{equation*}%
is a continuous embedding with dense range.

Further, if the unit ball of $B,$ 
\begin{equation*}
B_{\mathbf{1}}=\{x\in B:\left\Vert x\right\Vert _{B}\leq 1\},
\end{equation*}%
is closed in $A$ in projective topology, then 
\begin{equation*}
B=A_{b},
\end{equation*}%
where by $A_{b}$ we mean the bounded part of $A.$
\end{theorem}

\begin{proof}
Let $A$ be a real locally C*-algebra. We show that the bounded part $A_{b}$
of $A$\ is a real C*-algebra with required embedding. According to \cite%
{KatzFriedman06} 
\begin{equation*}
\mathfrak{A}=A\dotplus iA,
\end{equation*}%
is a complex locally C*-algebra, and there exists an involutive
antiautomorphism 
\begin{equation*}
\widehat{\Psi }:\mathfrak{A\rightarrow A},
\end{equation*}%
of order $2$ on $\mathfrak{A}$, so that 
\begin{equation*}
A=\{x\in \mathfrak{A}:\widehat{\Psi }(x)=\widehat{\Psi }(x^{\ast })\}.
\end{equation*}%
Let $\mathfrak{A}_{b}$ be the bounded part of $\mathfrak{A}$. Then,
according to \cite{Inoue71} $\mathfrak{A}_{b}$ is a complex C*-algebra with
a norm 
\begin{equation*}
\left\Vert .\right\Vert _{\infty }\equiv \left\Vert .\right\Vert _{\mathfrak{%
A}_{b}}.
\end{equation*}%
Let 
\begin{equation*}
\Psi :\mathfrak{A}_{b}\mathfrak{\rightarrow A}_{b},
\end{equation*}%
be a restriction of $\widehat{\Psi }$ to $\mathfrak{A}_{b}.$ Because $%
\mathfrak{A}_{b}$ is a *-subalgebra of $\mathfrak{A},$ $\Psi $ is as well an
involutive antiautomorphism of order $2$ on $\mathfrak{A}_{b},$ and, one can
see that the set 
\begin{equation*}
\{x\in \mathfrak{A}_{b}:\Psi (x)=\Psi (x^{\ast })\},
\end{equation*}%
is real isometrically *-isomorphic to the real C*-subalgebra of $\mathfrak{A}%
_{b}.$ It is a routine exercise to check that the bounded part $A_{b}$ of
the algebra $A$ is exactly equal to it: 
\begin{equation*}
A_{b}=\{x\in \mathfrak{A}_{b}:\Psi (x)=\Psi (x^{\ast })\}.
\end{equation*}

Conversely, let $A$ be a real lmc*-algebra with dense real *-subalgebra $B$
which is a real C*-algebra under a norm 
\begin{equation*}
\left\Vert .\right\Vert _{B}.
\end{equation*}%
Let 
\begin{equation*}
\mathfrak{A}=A\dotplus iA,
\end{equation*}%
and 
\begin{equation*}
\mathfrak{B}=B\dotplus iB.
\end{equation*}%
Similar to the arguments in \cite{KatzFriedman06} one can establish that $%
\mathfrak{A}$ is a complex lmc *-algebra. From the definition of real
C*-algebras (see \cite{Li03})\ it follows that $\mathfrak{B}$ is a complex
C*-algebra with the norm 
\begin{equation*}
\left\Vert .\right\Vert _{\mathfrak{B}},
\end{equation*}%
so that 
\begin{equation*}
\left\Vert a+ib\right\Vert _{\mathfrak{B}}=\sqrt{\left\Vert a\right\Vert
_{B}^{2}+\left\Vert b\right\Vert _{B}^{2}},
\end{equation*}%
$a,b\in B$. Also, from the construction of the complexification for the
algebra $A$ it follows that from the continuity of the embedding of $B$ into 
$A$ we can conclude a continuity of the embedding of $\mathfrak{B}$ into $%
\mathfrak{A}$. Thus, from Theorem 1 of Bhatt and Karia above (see as well 
\cite{BhattKaria93}) it follows that $\mathfrak{A}$ is a locally C*-algebra.
Because 
\begin{equation*}
\mathfrak{A}=A\dotplus iA,
\end{equation*}%
and 
\begin{equation*}
A\cap iA=\{\mathbf{0}\},
\end{equation*}%
from \cite{KatzFriedman06} it follows that $A$ is a real locally C*-algebra.

Finally, because the cussedness of 
\begin{equation*}
B_{\mathbf{1}}=\{x\in B:\left\Vert x\right\Vert _{B}\leq 1\},
\end{equation*}%
in the projective topology of $A$ is, due to the aforementioned
complexification arguments, equivalent to the cussedness of 
\begin{equation*}
\mathfrak{B}_{\mathbf{1}}=\{x\in \mathfrak{B}:\left\Vert x\right\Vert _{%
\mathfrak{B}}\leq 1\},
\end{equation*}%
in the projective topology of $\mathfrak{A},$ the proof is now completed.
\end{proof}

\section{Intrinsic characterization of locally JB-algebras}

In the current section we present a Jordan-algebraic version of Theorem 1.

\begin{theorem}
Let $M$ be a complete Jordan fine locally multiplicatively convex algebra.
Then $M$ is a locally JB-algebra iff $M$ contains a Jordan subalgebra $N$
such that:

1. $N$ is a JB-algebra with some norm $\left\Vert .\right\Vert _{N};$

\ \ \ \ \ \ \ \ \ \ \ \ \ \ \ \ \ \ \ \ \ \ \ \ \ \ \ \ \ \ \ \ \ \ \ \ \ \
\ \ \ \ and

2. the inclusion 
\begin{equation*}
(N,\left\Vert .\right\Vert _{N})\rightarrow M,
\end{equation*}%
is a continuous embedding with dense range.

Further, if the unit ball of $N,$ 
\begin{equation*}
N_{\mathbf{1}}=\{x\in N:\left\Vert x\right\Vert _{N}\leq 1\},
\end{equation*}%
is closed in $M$ in projective topology, then 
\begin{equation*}
N=M_{b},
\end{equation*}%
where by $M_{b}$ we mean the bounded part of $M.$
\end{theorem}

\begin{proof}
\qquad Let $M$ be a locally JB-algebra. Using functional calculus, one can
see that 
\begin{equation*}
x\circ (\mathbf{1}+x^{2})^{-1}\in M_{b},
\end{equation*}%
and 
\begin{equation*}
\left\Vert x\circ (\mathbf{1}+x^{2})^{-1}\right\Vert \leq 1.
\end{equation*}%
Then, for each $n\in 
%TCIMACRO{\U{2115} }%
%BeginExpansion
\mathbb{N}
%EndExpansion
,$ let us set 
\begin{equation*}
x_{n}=x\circ (\mathbf{1}+\frac{x^{2}}{n})^{-1}.
\end{equation*}%
It is easy to see that $x_{n}\in M_{b},$ for $\forall n\in 
%TCIMACRO{\U{2115} }%
%BeginExpansion
\mathbb{N}
%EndExpansion
.$ Now, for every $p$ be a continuous fine submultiplicative regular
seminorm on $M,$ 
\begin{equation*}
p(x-x_{n})\leq \frac{1}{\sqrt{n}}p(x)p(\frac{x}{\sqrt{n}}\circ (\mathbf{1}+%
\frac{x^{2}}{n})^{-1})\leq \frac{1}{\sqrt{n}}p(x)\rightarrow 0,
\end{equation*}%
as $n\rightarrow \infty ,$ which shows that $(M_{b},\left\Vert .\right\Vert
_{\infty })$ is continuously embedded in $M$ with sequentially dense range.

\qquad Conversely, let $M$ be a complete Jordan fine locally
multiplicatively convex algebra. Let 
\begin{equation*}
P=\{p_{\alpha }\},
\end{equation*}%
$\alpha \in \Lambda $ be a saturated separating directed family of
continuous fine submultiplicative seminorms on $M,$ and 
\begin{equation*}
M=\underset{\longleftarrow }{\lim }M_{\alpha },
\end{equation*}%
where 
\begin{equation*}
M_{\alpha }=M/J_{\alpha },
\end{equation*}%
with 
\begin{equation*}
J_{\alpha }=\{x\in M:p_{\alpha }(x)=0\},
\end{equation*}

$\alpha \in \Lambda .$

\ \ To begin with, let us assume that $M_{b}$ is a JB-algebra with a norm 
\begin{equation*}
\left\Vert x\right\Vert _{\infty }=\underset{\alpha \in \Lambda }{\sup }%
p_{\alpha }(x)<\infty ,
\end{equation*}%
dense in $M$ in projective topology of $M$. Let us replace the family $P$
with another saturated separating directed family 
\begin{equation*}
P^{^{\prime }}=\{p_{\alpha }^{^{\prime }}\},
\end{equation*}%
$\alpha \in \Lambda $ of continuous fine submultiplicative seminorms on $M,$
with 
\begin{equation*}
\underset{\alpha \in \Lambda }{\sup }p_{\alpha }(x)=\underset{\alpha \in
\Lambda }{\sup }p_{\alpha }^{^{\prime }}(x),
\end{equation*}%
\begin{equation*}
M=\underset{\longleftarrow }{\lim }M_{\alpha }^{^{\prime }},
\end{equation*}%
where 
\begin{equation*}
M_{\alpha }^{^{\prime }}=M/J_{\alpha }^{^{\prime }},
\end{equation*}%
with 
\begin{equation*}
J_{\alpha }^{^{\prime }}=\{x\in M:p_{\alpha }^{^{\prime }}(x)=0\},
\end{equation*}%
the following way: 
\begin{equation*}
p_{\alpha }^{^{\prime }}(x)=\left\{ 
\begin{array}{c}
\left\Vert x\right\Vert _{\infty },\text{ }\forall x\in M_{b}; \\ 
p_{\alpha }(x),\text{ }\forall x\notin M_{b};%
\end{array}%
\right. ,
\end{equation*}%
$\alpha \in \Lambda .$ One can easily see that with the norm 
\begin{equation*}
\left\Vert x_{\alpha }\right\Vert _{\alpha }=p_{\alpha }^{^{\prime }}(x),
\end{equation*}%
$M_{\alpha }^{^{\prime }}$ becomes a fine Banach-Jordan algebra, where 
\begin{equation*}
x_{\alpha }=x+J_{\alpha }^{^{\prime }}.
\end{equation*}

For a given $\alpha \in \Lambda ,$ let 
\begin{equation*}
\pi _{\alpha }:M\rightarrow M_{\alpha }^{^{\prime }}:
\end{equation*}%
\begin{equation*}
x\rightarrow x_{\alpha }=x+J_{\alpha }^{^{\prime }},
\end{equation*}%
be a continuous projection from $M$ onto $M_{\alpha }^{^{\prime }}.$ Because 
$p_{\alpha }^{^{\prime }}(x)$ is regular on $M_{b},$ $\pi _{\alpha }(M_{b})$
will be, due to projective topological density of $M_{b}$ in $M$, a dense in 
$\left\Vert .\right\Vert _{\alpha }$ norm Jordan subalgebra of $M_{\alpha
}^{^{\prime }},$ which is a pre-JB-algebra in the norm $\left\Vert
.\right\Vert _{\alpha }.$ It is remained to show that the norm $\left\Vert
.\right\Vert _{\alpha }$ is regular on $M_{\alpha }^{^{\prime }}$, i.e. that 
\begin{equation*}
\forall y_{\alpha }\in M_{\alpha }^{^{\prime }},\left\Vert y_{\alpha
}^{2}\right\Vert _{\alpha }=\left\Vert y_{\alpha }\right\Vert _{\alpha }^{2}.
\end{equation*}%
On the elements of $\pi _{\alpha }(M_{b})$, 
\begin{equation*}
\left\Vert x_{\alpha }^{2}\right\Vert _{\alpha }=\left\Vert x_{\alpha
}\right\Vert _{\alpha }^{2},
\end{equation*}%
and if 
\begin{equation*}
y_{\alpha }\notin \pi _{\alpha }(M_{b}),
\end{equation*}%
\begin{equation*}
y_{\alpha }\in M_{\alpha }^{^{\prime }},
\end{equation*}%
from submultiplicativity of $\left\Vert .\right\Vert _{\alpha }$ it already
follows that 
\begin{equation*}
\left\Vert y_{\alpha }^{2}\right\Vert _{\alpha }\leq \left\Vert y_{\alpha
}\right\Vert _{\alpha }^{2}.
\end{equation*}%
Let us assume on the contrary that 
\begin{equation*}
\left\Vert y_{\alpha }^{2}\right\Vert _{\alpha }-\left\Vert y_{\alpha
}\right\Vert _{\alpha }^{2}=\varepsilon >0.
\end{equation*}%
Then, because $\pi _{\alpha }(M_{b})$ is dense in $M_{\alpha }^{^{\prime }}$
in the norm $\left\Vert .\right\Vert _{\alpha },$ we can find 
\begin{equation*}
x_{\alpha }\in \pi _{\alpha }(M_{b}),
\end{equation*}%
to be such that 
\begin{equation*}
\left\Vert y_{\alpha }^{2}\right\Vert _{\alpha }-\left\Vert x_{\alpha
}^{2}\right\Vert _{\alpha }<\frac{\varepsilon }{2},
\end{equation*}%
and 
\begin{equation*}
\left\Vert y_{\alpha }\right\Vert _{\alpha }^{2}-\left\Vert x_{\alpha
}\right\Vert _{\alpha }^{2}<\frac{\varepsilon }{2}.
\end{equation*}%
Thus, we have 
\begin{equation*}
\varepsilon =\left\Vert y_{\alpha }^{2}\right\Vert _{\alpha }-\left\Vert
y_{\alpha }\right\Vert _{\alpha }^{2}=\left\Vert y_{\alpha }^{2}\right\Vert
_{\alpha }-\left\Vert x_{\alpha }^{2}\right\Vert _{\alpha }+\left\Vert
x_{\alpha }\right\Vert _{\alpha }^{2}-\left\Vert y_{\alpha }\right\Vert
_{\alpha }^{2}=
\end{equation*}

\begin{equation*}
=|\left\Vert y_{\alpha }^{2}\right\Vert _{\alpha }-\left\Vert x_{\alpha
}^{2}\right\Vert _{\alpha }|+|\left\Vert x_{\alpha }\right\Vert _{\alpha
}^{2}-\left\Vert y_{\alpha }\right\Vert _{\alpha }^{2}|<\varepsilon .
\end{equation*}%
Contradiction.

\ \ Now, let us assume that $N$ is a JB-algebra with some norm $\left\Vert
.\right\Vert _{N}$, and the inclusion 
\begin{equation*}
(N,\left\Vert .\right\Vert _{N})\rightarrow M,
\end{equation*}%
is a continuous embedding with dense in the projective topology $\tau $ of $%
M $\ range. Let $B_{M}(\tau )$ be the collection of all absolutely-convex,
closed, bounded, idempotent subsets of $M,$ containing identity $\mathbf{1}$
of $M$.

Let 
\begin{equation*}
M(S)=\{\lambda x:\lambda \in 
%TCIMACRO{\U{211d} }%
%BeginExpansion
\mathbb{R}
%EndExpansion
,x\in S\},
\end{equation*}%
and the Minkowski functional of $S$ in $M(S)$ is 
\begin{equation*}
\left\vert x\right\vert _{S}=\inf \{\lambda >0:x\in M(S)\}.
\end{equation*}%
If the set $S$ is absorbing in $M(S),$ one can see that 
\begin{equation*}
M(S)=M,
\end{equation*}%
and $\left\vert .\right\vert _{S}$ is defined on the whole of $M.$ Because $%
(M,\tau )$ is complete, then for 
\begin{equation*}
B\in B_{M}(\tau ),
\end{equation*}%
$(M(B),\left\vert x\right\vert _{B})$ is a Banach-Jordan algebra. One can
see that given a saturated separating directed family of continuous fine
submultiplicative seminorms on $M:$ 
\begin{equation*}
P=\{p_{\alpha }\},
\end{equation*}%
$\alpha \in \Lambda $, 
\begin{equation*}
U_{P}=\{x\in M:p_{\alpha }(x)\leq 1,\forall \alpha \in \Lambda \}\in
B_{M}(\tau ),
\end{equation*}%
and conversely, given a subset 
\begin{equation*}
B\in B_{M}(\tau ),
\end{equation*}%
there exists such a saturated separating directed family of continuous fine
submultiplicative seminorms on $M:$ 
\begin{equation*}
P=\{p_{\alpha }\},
\end{equation*}%
$\alpha \in \Lambda $, so that 
\begin{equation*}
B\subset U_{P}.
\end{equation*}

Now, it is clear that the closure $\overline{N_{\mathbf{1}}}$ of $N_{\mathbf{%
1}}$ in $(M,\tau )$ belongs to $B_{M}(\tau ),$ therefore there exists a
saturated separating directed family of continuous fine submultiplicative
seminorms on $M:$ 
\begin{equation*}
P=\{p_{\alpha }\},
\end{equation*}%
$\alpha \in \Lambda $, such that 
\begin{equation*}
N_{\mathbf{1}}\subset U_{P},
\end{equation*}%
so, 
\begin{equation*}
(M(U_{P}),\left\vert x\right\vert _{U_{P}}),
\end{equation*}%
is a Banach-Jordan algebra, and 
\begin{equation*}
N\subset (M(U_{P}),\left\vert x\right\vert _{U_{P}}),
\end{equation*}%
thus 
\begin{equation*}
\left\vert x\right\vert _{U_{P}}\leq \left\Vert x\right\Vert _{N}.
\end{equation*}%
Moreover, because $(N,\left\Vert .\right\Vert _{N})$ is a JB-algebra, using
spectral radii of an element in $N$ and $M(U_{P})$, one can conclude that 
\begin{equation*}
\left\Vert x\right\Vert _{N}=\left\vert x\right\vert _{U_{P}},
\end{equation*}%
for all $x\in N.$

Due to the fact that $(N,\left\Vert .\right\Vert _{N})$ is dense in $(M,\tau
),$ the same way as above we now can establish that $M$ is a locally
JB-algebra.

Finally, if $N_{\mathbf{1}}$ is closed in $(M,\tau ),$ then 
\begin{equation*}
(M(U_{P}),\left\vert x\right\vert _{U_{P}})=(N,\left\Vert .\right\Vert _{N}).
\end{equation*}%
Indeed, if 
\begin{equation*}
x\in M(U_{P})\subset M,
\end{equation*}%
then 
\begin{equation*}
x_{n}=x\circ (\mathbf{1}+\frac{x^{2}}{n})^{-1}\in N\subset M(U_{P}),
\end{equation*}%
and 
\begin{equation*}
\left\vert x-x_{n}\right\vert _{U_{P}}\leq \frac{1}{\sqrt{n}}\left\vert
x\right\vert _{U_{P}}\rightarrow 0,
\end{equation*}%
as $n\rightarrow \infty ,$ therefore $x\in N.$ Thus 
\begin{equation*}
M(U_{P})=N=M_{b}.
\end{equation*}%
The proof is now complete.
\end{proof}


\begin{thebibliography}{99}
\bibitem{AlfsenSchultzStoermer78} \textbf{Alfsen, E.M.; Shultz, F.W.; St\o %
rmer, E.,} \textit{A Gelfand-Naimark theorem for Jordan algebras.} (English)
Advances in Math. Vol. 28 (1978), No. 1, pp. 11-56.

\bibitem{Apostol71} \textbf{Apostol, C.,} \textit{b*-algebras and their
representation.} (English) J. London Math. Soc. (2) No. 3 (1971), pp. 30--38.

\bibitem{Arens52} \textbf{Arens, R.,} \textit{A generalization of normed
rings.} (English), Pac. J. Math., Vol. 2 (1952), pp. 455-471.

\bibitem{BhattKaria93} \textbf{Bhatt, S.J.; Karia, D.J.,} \textit{An
intrinsic characterization of pro-C*-algebras and its applications.}
(English summary) J. Math. Anal. Appl. Vol. 175 (1993), No. 1, pp. 68--80.

\bibitem{Fragoulopoulou05} \textbf{Fragoulopoulou, M.,} \textit{Topological
algebras with involution.} (English) North-Holland Mathematics Studies, Vol.
200. Elsevier Science B.V., Amsterdam, 495 pp., (2005).

\bibitem{GilesKoehler73} \textbf{Giles, J. R.; Koehler, D. O.,} \textit{On
numerical ranges of elements of locally m-convex algebras.} (English)
Pacific J. Math. 49 (1973), pp. 79--91.

\bibitem{GelfandNaimark43} \textbf{Gelfand, I.M.; Naimark, M.A.,} \textit{On
the embedding of normed rings into the ring of operators in Hilbert space.}
(English. Russian summary) Rec. Math. [Mat. Sbornik] N.S. Vol. 12(54)
(1943), pp. 197-213.

\bibitem{Hanche-OlsenStoermer84} \textbf{Hanche-Olsen, H.; St\o rmer, E.,} 
\textit{Jordan operator algebras.} (English), Monographs and Studies in
Mathematics, Vol. 21. Boston - London - Melbourne: Pitman Advanced
Publishing Program. VIII, 183 pp., (1984).

\bibitem{Inoue71} \textbf{Inoue, A.,} \textit{Locally C*-algebras.}
(English), Mem. Fac. Sci. Kyushu Univ. (Ser. A), No. 25 (1971), pp. 197-235.

\bibitem{KatzFriedman06} \textbf{Katz, A.A.; Friedman, O.,} \textit{On
projective limits of real C*- and Jordan operator algebras.} (English),
Vladikavkaz Mathematical Journal, Vol. 8 (2006), No. 2, pp. 33-38.

\bibitem{Li03} \textbf{Li, B.,} \textit{Real operator algebras.} (English)
World Scientific Publishing Co., Inc., River Edge, NJ, 241 pp., (2003).

\bibitem{Michael52} \textbf{Michael, E.A.,} \textit{Locally
multiplicatively-convex topological algebras.} (English), Mem. Am. Math.
Soc., Vol. 11 (1952), 79 pp.

\bibitem{Pedersen79} \textbf{Pedersen, G.K.,} \textit{C*-algebras and their
automorphism groups.} (English), London Mathematical Society Monographs.
Vol. 14.London - New York -San Francisco: Academic Press., 416 pp., (1979).

\bibitem{Phillips88} \textbf{Phillips, N. C.,} \textit{Inverse limits of
C*-algebras.} (English) J. Operator Theory Vol. 19 (1988), No. 1, pp.
159--195.

\bibitem{Schmuedgen75} \textbf{Schm\"{u}dgen, K.,} \textit{\"{U}ber
LMC-Algebren.} (German) Math. Nachr.Vol. 68 (1975), pp. 167--182.
\end{thebibliography}
\end{document}